\documentclass{siamart190516}
\usepackage{hyphenat}

\let\theoremstyle\relax

\usepackage{amsmath,amsfonts,amsthm, amssymb, bm}
\usepackage{array}
\newcolumntype{P}[1]{>{\centering\arraybackslash}p{#1}}
\title{Perturbing ordinary differential equations to generate resonant and repeated root solutions\thanks{Submitted to the editors January 31, 2021.
\funding{B.G. is supported by the Paul and Daisy Soros Fellowship and the NSF Graduate Research Fellowship Program.}}}
\author{Bernardo Gouveia\thanks{Department of Chemical and Biological Engineering, Princeton University, Princeton, NJ 08544, USA.}
\and Howard A. Stone\thanks{Department of Mechanical and Aerospace Engineering, Princeton University, Princeton, NJ 08544, USA. (\email{hastone@princeton.edu}).}}

\headers{ODEs: resonance and repeated roots}{Bernardo Gouveia and Howard A. Stone}

\begin{document}
\maketitle

\begin{abstract}
In the study of ordinary differential equations (ODEs) of the form $\hat{L}[y(x)]=f(x)$, where $\hat{L}$ is a linear differential operator, two related phenomena can arise: resonance, where $f(x)\propto u(x)$ and $\hat{L}[u(x)]=0$, and repeated roots, where $f(x)=0$ and $\hat{L}=\hat{D}^n$ for $n\geq 2$. We illustrate a method to generate exact solutions to these problems by taking a known homogeneous solution $u(x)$, introducing a parameter $\epsilon$ such that $u(x)\rightarrow u(x;\epsilon)$, and Taylor expanding $u(x;\epsilon)$ about $\epsilon = 0$. The coefficients of this expansion $\frac{\partial^k u}{\partial\epsilon^k}\big{|}_{\epsilon=0}$  yield the desired resonant or repeated root solutions to the ODE. This approach, whenever it can be applied, is more insightful and less tedious than standard methods such as reduction of order or variation of parameters. While the ideas can be introduced at the undergraduate level, we could not find any elementary or advanced text that illustrates these ideas with appropriate generality.
\end{abstract}



\section{Introduction}

Introductory courses on ODEs often focus on illustrating me\hyp{}thods of solution for problems in the natural sciences and engineering. One common topic is resonance, which for the typical case of a mechanical or electrical oscillator corresponds to forcing the system at its natural frequency. The unwanted oscillations that occurred when the Millennium Bridge, which crosses over the Thames River in London, was opened to foot traffic in June 2000 serves a modern realization \cite{strogatz2005crowd}. Mathematically, this response occurs because the governing ODE is forced with one of its homogeneous solutions. 

Of course, modern computer algebra systems such as \textit{Mathematica} and \textit{Maple} provide immediate, often analytical, solutions to these problems. However, because standard algorithms are utilized, it is not uncommon that the results are cumbersome and not insightful. For the special case of resonance, we show here using ideas grounded in perturbation theory and analytic continuation that the resonant solution can be obtained simply by constructing an appropriate Taylor series. This method yields simpler, more transparent functional forms for the resonant solution that have an obvious relationship to the homogeneous solution. We will also see that the 
same ideas are applicable when constructing linearly independent solutions of equations with ``repeated roots''. 
 
To illustrate the method, we first provide some general background on linear ODEs to orient the reader on the scope of the problems we seek to solve. We then proceed to go through a series of examples that explicitly demonstrate how the method is used. To conclude, we provide a general derivation that summarizes the method and reveals its underlying structure. Because this approach seems particularly flexible, practically requires only knowledge of Taylor series, and in non-elementary problems produces much simpler solutions than mathematical software packages, we believe it will be of interest. These ideas can be taught in a course at the undergraduate level, although we are not aware of any book on ODEs at \textit{any} level that illustrates the approach nor emphasizes its generality. 

\section{Background}

In the study of $n^{\mathrm{th}}$-order linear ODEs of the form 
\begin{equation}
\hat{L}\left[y(x)\right]=f(x), 
\quad\hbox{where}\quad \hat{L}=a_0(x) + \sum_{j=1}^{n}a_j(x)\frac{\mathrm{d}^j}{\mathrm{d}x^j},
\end{equation}solutions may be represented as $y(x) = \sum_{j=1}^{n}c_ju_j(x) + u_p(x)$. Here, $u_j(x)$ are the $n$ linearly independent homogeneous solutions such that  $\hat{L}\left[u_j(x)\right] = 0$ and $u_p(x)$ is a particular solution such that $\hat{L}\left[u_p(x)\right] = f(x)$. 
The constants $\lbrace c_j\rbrace$ are determined by auxiliary data. 
Two common special cases arise when seeking the solution set $\lbrace u_j(x), u_p(x) \rbrace$: resonance and repeated roots.  

The first, resonance, occurs when $f(x) \propto u_j(x)$, i.e., when the forcing function includes any term proportional to one of the homogeneous solutions. Therefore, positing a particular solution of the form $u_p(x)\propto f(x)$, an intuitive ansatz for a linear system, will automatically fail since $\hat{L}\left[f(x)\right]=0$. 

The second is the case of repeated roots, which concerns the homogeneous problem $f(x)=0$ whenever the linear operator can be factored as $\hat{L}=\prod_{k=1}^{m<n}\hat{D}_k$, where all of the $\hat{D}_k$ commute. In this case it is necessary
to solve the lower order sub-problems $\hat{D}_k\left[y(x)\right] = 0$. The difficulty arises when one of the operators, say $\hat{D}_1$, is repeated $r$ times, and so solving the sub-problems can produce only $n-r-1$ linearly independent solutions. 

Both these cases are systematically resolved by reduction of order\footnote{Variation of parameters is the most general method to solve all nonhomogeneous linear ODEs, but it is just a logical extension of reduction of order.}, which is a standard topic in introductory courses \cite[sec. 3.4, 3.8]{boyce2017}, \cite[sec. 2.8, 2.10]{Kreyszig2005}, \cite[sec. 6.7]{Jeffrey2001}, \cite[sec. 3.4.4, 3.6.2]{greenberg1988}. Substituting a solution of the form $y(x) = u(x)g(x)$, where $u(x)$ is an already obtained homogeneous solution, into the ODE results in a reduced-order ODE for $g(x)$ that is easier to solve and will generate the remaining resonant or repeated root solutions. In practice, the challenge with applying this method is that it only reduces the order of the ODE by one, and so it typically fails to be useful for ODEs that are higher than second order. Also, for non-elementary ODEs, the method may only guarantee complicated integral representations, which are difficult to put into a simpler form. Furthermore, even when the method is successful, it is a tedious calculation and the connection between the resonant or repeated root solution and the known homogeneous solution $u_j(x)$ is not obvious. 

Here, we highlight what we believe is a better, or at least more insightful, way to obtain resonant or repeated root solutions to linear ODEs. For the case of resonance, the idea is to introduce a small parameter $\epsilon$ in the forcing function so that $f(x) \rightarrow f(x; \epsilon)$ such that $\lim_{\epsilon\rightarrow 0}f(x; \epsilon) = f(x)$. It will then be possible to construct the particular solution by positing $u_p(x; \epsilon) \propto f(x; \epsilon)$. Taylor expanding $u_p(x; \epsilon)$ in $\epsilon$ and setting $\epsilon=0$ will generate the desired resonant solution. In the examples, we will show that repeated root solutions may be constructed in the exact same way. In the end, we will learn that the obtained resonant or repeated root solution is simply an analytic continuation of the known homogeneous solution in $\epsilon$. This description is formal, but the method will become concrete in the following step-by-step examples.

The rudiments of this idea have been discussed in advanced texts \cite[pgs. 11-12]{bender2013advanced}, \cite[pgs. 17-18]{white2010asymptotic}, but only in the context of repeated roots for constant coefficient or equidimensional ODEs, where reduction to an algebraic characteristic equation is possible. As we will see, the method as presented here is more general and useful in a variety of non-elementary problems, where no algebraic characteristic equation exists globally. Makarov et. al \cite{makarov1978construction} give a more general presentation of the idea, so it is certainly known in the Russian mathematical literature, but even there the discussion is limited to the case of resonance with no mention of repeated root solutions. 

While our method only applies to problems involving resonance or repeated roots, such problems occur often enough in applied mathematics and the physical sciences that we believe it is worth highlighting. Furthermore, as previously noted, this approach is not presented at all in introductory texts on differential equations to the best of our knowledge, but is certainly comprehensible at the undergraduate level and is just as useful as reduction of order. Whenever one encounters resonance or a repeated root, we suggest that the method we now illustrate be the method of choice.

\section{Examples}

In this section we give seven concrete examples, in order of increasing complexity. We suggest reading them in order, as the ideas naturally build off each other. The first example of resonance in a constant coefficient $2^{\mathrm{nd}}$-order ODE is discussed in many introductory texts using reduction of order \cite{boyce2017, Kreyszig2005, Jeffrey2001, greenberg1988}. 

\subsection{Constant coefficient equation at resonance} \label{example1}

Consider the constant coefficient ODE
\begin{equation} \label{ode:1}
\frac{\mathrm{d}^2y}{\mathrm{d}x^2} + y = \sin x.
\end{equation}
The homogeneous solutions are $u_1(x) = \sin x$ and $u_2 (x) = \cos x$. We observe that the forcing function $f(x) = \sin x$ is linearly dependent on $u_1(x)$, in fact \textit{it is} $u_1(x)$, and so we have resonance. Therefore, naively guessing $u_p(x)=A\sin x$, where $A$ is a to be determined coefficient, will fail as it would render the left-hand side of equation (\ref{ode:1}) equal to zero. The strategy is to then introduce a parameter $\epsilon$ to equation (\ref{ode:1}) such that
\begin{equation} \label{ode:1_parameter}
\frac{\mathrm{d}^2y}{\mathrm{d}x^2} + y = \sin\left((1+\epsilon)x\right),
\end{equation}
where we are interested in the limit  $\epsilon\rightarrow0$. Now, substituting $u_p(x)=A\sin\left((1+\epsilon)x\right)$ provides a particular solution so long as $A=\frac{1}{1 - (1+\epsilon)^2}$. Hence, the general solution is
\begin{equation} \label{ode:1_solution}
y(x) = c_1\sin x + c_2\cos x - \frac{\sin\left(x + \epsilon x\right)}{2\epsilon + \epsilon^2}.
\end{equation}
Constructing a Taylor expansion around $\epsilon=0$, and neglecting terms of $\mathrm{O}\left(\epsilon^2\right)$ or higher, gives 
\begin{equation} \label{ode:1_solution_continued}
y(x) = \left(c_1 -\frac{1}{2\epsilon}\right)\sin x+ c_2\cos x - \frac{x\cos x}{2}.
\end{equation}
Because we have complete freedom in choosing $c_1$, we can make the transformation $c_1 = \tilde{c}_1 - 1/2\epsilon$, which removes the divergence as $\epsilon\rightarrow 0$.  Thus, we have obtained the general solution 
\begin{equation} \label{ode:1_final_solution}
y(x) = \tilde{c}_1\sin x + c_2\cos x - \frac{x\cos x}{2}.
\end{equation}
If $x$ represents time, the $x\cos x$ behavior is the signature of the ever growing oscillations of a conservative oscillator at resonance. 

Let us take stock of the solution approach. The first step is to introduce a parameter $\epsilon$ that allows us to guess a particular solution $u_p(x; \epsilon) \propto f(x; \epsilon)$, where the limit $\epsilon\rightarrow0$ corresponds to the problem we wish to solve. This limit will be singular, but Taylor expanding $u_p(x; \epsilon)$ to linear order in $\epsilon$ allows grouping the singular part with a homogeneous solution. Relabeling a free constant removes the singularity and gives the desired $\epsilon=0$ solution. 

\subsection{Repeated roots of the equidimensional equation} \label{example2}

Consider the equi\hyp{}dimensional ODE in the form
\begin{equation} \label{ode:2}
x^2\frac{\mathrm{d}^2y}{\mathrm{d}x^2} + \left(1 - 2b\right)x\frac{\mathrm{d}y}{\mathrm{d}x} + b^2y = 0,
\end{equation}
where $b$ is a given constant. Carefully staring at this equation, we observe that the linear operator can be factored and the ODE can be rewritten as
\begin{equation} \label{ode:2_factored}
\left(x\frac{\mathrm{d}}{\mathrm{d}x} - b\right)^2y = 0,
\end{equation}
and thus this equation has a repeated root. We construct the homogeneous solutions $\left \{ u_1(x), u_2(x)\right \}$ by solving the sub-problems $\left(x\frac{\mathrm{d}}{\mathrm{d}x} - b\right)u_1=0$ and $\left(x\frac{\mathrm{d}}{\mathrm{d}x} - b\right)u_2=u_1$. By inspection we observe that both $u_1(x)$ and $u_2(x)$ are indeed solutions to equation (\ref{ode:2_factored}). Integration gives $u_1(x) = x^b$ and we are left to solve
\begin{equation}  \label{ode:2_subproblem}
\left(x\frac{\mathrm{d}}{\mathrm{d}x} - b\right)u_2 = x^b.
\end{equation}

We are now in the same position as example \hyperref[example1]{3.1}, since the homogeneous solution of equation (\ref{ode:2_subproblem}) is exactly the forcing term $x^b$, and so we have resonance. Thus, we learn that resonance and repeated roots are effectively the same feature. Writing $b\rightarrow b+\epsilon$ gives
\begin{equation}  \label{ode:2_parameter}
\left(x\frac{\mathrm{d}}{\mathrm{d}x} - b\right)u_2 = x^{b+\epsilon}.
\end{equation}The equation and solution of interest correspond to the limit $\epsilon\rightarrow 0$. 
Substituting $u_2(x) = Ax^{b+\epsilon}$ as the particular solution is successful if $A=1/\epsilon$, and therefore the solution to equation (\ref{ode:2_factored}) can be written
\begin{equation} \label{ode:2_solution_continued}
y(x) = c_1x^b + c_2\frac{x^{b}x^{\epsilon}}{\epsilon}.
\end{equation}
Constructing the Taylor expansion for $0<\epsilon \ll 1$, we have $x^\epsilon = 1 + \epsilon\log x + \mathrm{O}\left(\epsilon^2\right)$, therefore
\begin{equation} \label{ode:2_solution_singular}
y(x) = \left(c_1 + \frac{c_2}{\epsilon}\right)x^b + c_2x^b\log x.
\end{equation}
Relabeling the constant $c_1 = \tilde{c}_1 - c_2/\epsilon$ removes the divergence and gives the desired solution
\begin{equation} \label{ode:2_final_solution}
y(x) = \tilde{c}_1x^b + c_2x^b\log x.
\end{equation}

Of course, in this particular example, equation (\ref{ode:2_subproblem})  could have been solved more directly via an integrating factor. However, that approach only works for first-order ODEs, whereas our method easily generalizes to higher-order ODEs, as illustrated in example \hyperref[example4]{3.4}. 

\subsection{Airy's equation at resonance} \label{example3}

So far, we have only looked at ODEs that may be reduced to algebraic characteristic equations, the solutions to which yield elementary functions. Furthermore, we have yet to consider a well-posed initial (IVP) or boundary (BVP) value problem. To increase the complexity and demonstrate the generality of our method, consider the forced Airy BVP
\begin{subequations} \label{ode:3-airy}
	\begin{eqnarray}
	\frac{1}{x}\frac{\mathrm{d}^2y}{\mathrm{d}x^2} - y &=& \mathrm{Ai}(x) \\
	y(0) &=& 1 \\
	y(x\rightarrow\infty)&\rightarrow& 0,
	\end{eqnarray}
\end{subequations}where $\mathrm{Ai}(x)$ is the Airy function of the first kind. 
Similar equations come up when studying the Schr\"{o}dinger equation for a particle in a linear potential \cite{griffiths2018introduction}. The homogeneous solutions are the two linearly independent Airy functions of the first and second kind, respectively,  $u_1(x)=\mathrm{Ai}(x)$ and $u_2(x)=\mathrm{Bi}(x)$, and thus equation (\hyperref[ode:3-airy]{3.13a}) corresponds to resonance since the forcing function is proportional to a homogeneous solution. 

Following the same steps as in the first two examples, we introduce a parameter $\epsilon$ so that
\begin{equation} \label{ode:3_continued}
\frac{1}{x}\frac{\mathrm{d}^2y}{\mathrm{d}x^2} - y = \mathrm{Ai}\left((1+\epsilon)x\right),
\end{equation}
where we will let $\epsilon\rightarrow0$. Proposing a particular solution $u_p (x) = C\mathrm{Ai}\left((1+\epsilon)x\right)$ is successful if $C=\frac{1}{(1+\epsilon)^3-1}$, where we have used $\frac{1}{x}\frac{\mathrm{d}^2}{\mathrm{d}x^2}\left[\mathrm{Ai}(cx)\right]=c^3\mathrm{Ai}(cx)$. Therefore, the solution to (\ref{ode:3_continued}) is
\begin{equation} \label{ode:3_continued_solution}
y(x) = c_1\mathrm{Ai}(x) + c_2\mathrm{Bi}(x) + \frac{\mathrm{Ai}\left((1+\epsilon)x\right)}{\left(1+\epsilon\right)^3-1}.
\end{equation}
Constructing the Taylor expansion $\mathrm{Ai}(x + \epsilon x)=\mathrm{Ai}(x) + \epsilon x\mathrm{Ai}'(x) + \mathrm{O}\left(\epsilon^2\right)$, where $\frac{\mathrm{d}}{\mathrm{d}x}=~'$, and retaining terms only up to first order in $\epsilon$ gives
\begin{equation} \label{ode:3_continued_solution_singular}
y (x) = \left(c_1  + \frac{1}{3\epsilon}\right)\mathrm{Ai}(x) + c_2\mathrm{Bi}(x) + \frac{x\mathrm{Ai}'(x)}{3}.
\end{equation}
We can now apply the boundary data. The condition (\hyperref[ode:3-airy]{3.13c}) forces $c_2=0$ since $\mathrm{Bi}(x\rightarrow\infty)\rightarrow\infty$. The condition (\hyperref[ode:3-airy]{3.13b}) results in $1 = \left(c_1+\frac{1}{3\epsilon}\right)\mathrm{Ai}(0)\implies c_1 + \frac{1}{3\epsilon}= \frac{2\pi 3^{1/6}}{\Gamma(1/3)}$, where we have used $\mathrm{Ai}(0)=\frac{\Gamma(1/3)}{2\pi 3^{1/6}}$ (see \hyperref[Apdx-1]{Appendix 5.1} for a derivation). Thus we see that the divergence at $\epsilon=0$ goes away naturally when we consider a well-posed BVP\footnote{We thank Dionisios Margetis for emphasizing this point to us.}, and the final solution is
\begin{equation} \label{ode:3_solution}
y(x) = \frac{2\pi 3^{1/6}}{\Gamma(1/3)}\mathrm{Ai}(x) + \frac{x\mathrm{Ai}'(x)}{3}.
\end{equation}

For problems involving special functions, it is often useful to utilize \textit{Mathematica} as a first check on solvability. Entering the appropriate commands results in
{\small
\begin{subequations} \label{ode:3_DSolve}
	\begin{eqnarray}
	\texttt{sol} &=& \texttt{DSolve[(1/x)*y''[x] - y[x] == AiryAi[x], y[x], x]} \nonumber\\
	\texttt{y[x]} &\rightarrow& \texttt{1/6(}\pi\texttt{AiryAi[x]AiryAiPrime[x]AiryBi[x] - 2}\pi\texttt{xAiryAiPrime[x]}^2\nonumber \\
	\texttt{AiryBi[x]}&\texttt{-}&\pi\texttt{AiryAi[x]}^2\texttt{AiryBiPrime[x] + 2}\pi\texttt{xAiryAi[x]AiryAiPrime[x]} \nonumber\\
	\texttt{AiryBiPrime[x])}&+&\texttt{AiryAi[x]c[1] + AiryBi[x]c[2]} \nonumber,
	\end{eqnarray}
\end{subequations}
}which is rather disastrous compared to our solution! Because the structure of \textit{Mathematica}'s solution involves products of $\mathrm{Ai}(x)$ and $\mathrm{Bi}(x)$, it is clear that reduction of order is the algorithm underpinning this solution. Yet, we are guaranteed that the particular solution is unique up to added multiples of the homogeneous solutions (the Fredholm alternative). Indeed, when we apply \textit{Mathematica}'s brute-force simplification algorithm to this solution we find
{\small
\begin{subequations} \label{ode:3_FullSimplify}
	\begin{eqnarray}
	~~~~~\texttt{sol2} &=& \texttt{FullSimplify[sol]} \nonumber \\
	\texttt{y[x]} &\rightarrow& \texttt{1/3xAiryAiPrime[x] + AiryAi[x](-1/6 + c[1]) + AiryBi[x]c[2]} \nonumber,
	\end{eqnarray}
\end{subequations}
}which is precisely our general solution (\ref{ode:3_continued_solution_singular}) with a relabeling of a constant. 

Using our method, one can now appreciate how to obtain the simplified solution directly from the ODE. There is no need to generate the complicated solution via reduction of order and use esoteric properties of Airy functions to simplify it further. 

\subsection{Repeated roots of a fourth-order Bessel-like equation} \label{example4}

Our method easily generalizes to higher-order ODEs. Consider the fourth-order equation
\begin{equation} \label{ode:4_bessel}
\left(\frac{\mathrm{d}^2}{\mathrm{d}x^2} - \frac{1}{x}\frac{\mathrm{d}}{\mathrm{d}x} - k^2\right)^2y = 0,
\end{equation}
which arises in the study of hydrodynamic stability \cite{goren1962instability}. Here $k$ is a given constant. Equation (\ref{ode:4_bessel})  has a repeated root, so we construct the homogeneous solutions by solving the sub-problems $\left(\frac{\mathrm{d}^2}{\mathrm{d}x^2} - \frac{1}{x}\frac{\mathrm{d}}{\mathrm{d}x} - k^2\right)u_1=0$ and $\left(\frac{\mathrm{d}^2}{\mathrm{d}x^2} - \frac{1}{x}\frac{\mathrm{d}}{\mathrm{d}x} - k^2\right)u_2=u_1$. Substituting $u_1(x)=x\zeta(x)$ into the first sub-problem furnishes $x^2\frac{\mathrm{d}^2\zeta}{\mathrm{d}x^2} + x\frac{\mathrm{d}\zeta}{\mathrm{d}x} - \left(1 + k^2x^2\right)\zeta = 0$, for which the solutions are the modified Bessel functions $\zeta(x) = c_1I_1(kx) + c_2K_1(kx)$. Therefore, $u_1 (x) = c_1xI_1(kx) + c_2xK_1(kx)$ and the second sub-problem becomes
\begin{equation} \label{ode:4_second_problem}
\left(\frac{\mathrm{d}^2}{\mathrm{d}x^2} - \frac{1}{x}\frac{\mathrm{d}}{\mathrm{d}x} - k^2\right)u_2=c_1xI_1(kx) + c_2xK_1(kx).
\end{equation}

Equation (\ref{ode:4_second_problem}) displays resonance, as the forcing functions are the homogeneous solutions of the differential operator. Therefore, we introduce $\epsilon$ so that
\begin{equation} \label{ode:4_second_problem_continued}
\left(\frac{\mathrm{d}^2}{\mathrm{d}x^2} - \frac{1}{x}\frac{\mathrm{d}}{\mathrm{d}x} - k^2\right)u_2=c_1xI_1\left((k+\epsilon)x\right) + c_2xK_1\left((k+\epsilon)x\right)
\end{equation}
and guess the particular solution $u_p(x)=AxI_1\left((k+\epsilon)x\right) + BxK_1\left((k+\epsilon)x\right)$. Substituting this guess works only if we choose $A=\frac{c_1}{(k+\epsilon)^2-k^2}$ and $B=\frac{c_2}{(k+\epsilon)^2-k^2}$. Therefore, the solution to equation (\ref{ode:4_bessel})  is
\begin{equation} \label{ode:4_continued_solution}
y(x) = c_1xI_1(kx) + c_2xK_1(kx) + c_3\frac{xI_1\left((k+\epsilon)x\right)}{(k+\epsilon)^2-k^2} + c_4\frac{xK_1\left((k+\epsilon)x\right)}{(k+\epsilon)^2-k^2},
\end{equation}
where we have renamed some integration constants. We now construct the Taylor expansions
\begin{subequations} \label{ode:4_TaylorSeriesBesselFunctions}
\begin{eqnarray}
I_1(kx + \epsilon x) &=& I_1(kx) + \epsilon\left(I_1(kx)/k + xI_2(kx)\right) + \mathrm{O}\left(\epsilon^2\right) \\
\hbox{and}~~K_1(kx + \epsilon x) &=& K_1(kx) + \epsilon\left(K_1(kx)/k - xK_2(kx)\right) + \mathrm{O}\left(\epsilon^2\right),
\end{eqnarray}
\end{subequations}
where we have made use of the formulas $I_1'(z) = I_1(z)/z + I_2(z)$ and $K_1'(z) = K_1(z)/z - K_2(z)$. Combining (\ref{ode:4_continued_solution}-\ref{ode:4_TaylorSeriesBesselFunctions}) and keeping only terms to linear order in $\epsilon$ gives, after some rearrangement,
\begin{subequations} \label{ode:4_linear_order}
    \begin{eqnarray}
    y(x) &=& \left(c_1 + \frac{c_3}{2k\epsilon} + \frac{c_3}{2k^2}\right)xI_1(kx) + \left(c_2 + \frac{c_4}{2k\epsilon} + \frac{c_4}{2k^2}\right)xK_1(kx) \\
    &+&\frac{c_3}{2k}x^2I_2(kx) - \frac{c_4}{2k}x^2K_2(kx) \nonumber.
    \end{eqnarray}
\end{subequations}
Relabeling all constants to remove the $\epsilon\rightarrow0$ divergence and to clean up the final solution results in 
\begin{equation} \label{ode:4_solution}
y(x) = \tilde{c}_1xI_1(kx) + \tilde{c}_2xK_1(kx) + \tilde{c}_3x^2I_2(kx) + \tilde{c_4}x^2K_2(kx).
\end{equation}

Calling \textit{Mathematica} again for comparison, we find
{\small
\begin{subequations} \label{ode:4_DSolve}
	\begin{eqnarray}
    \texttt{inner} &=& \texttt{D[y[x], \{x, 2\}] - (1/x)*D[y[x], x] - k}^2\texttt{*y[x]} \nonumber \\
	\texttt{sol} &=& \texttt{DSolve[D[inner, \{x, 2\}] - (1/x)*D[inner, x] - k}^2\texttt{*inner == 0, y[x], x]} \nonumber\\
	\texttt{y[x]} &\rightarrow& \texttt{x}^2\texttt{BesselJ[2, ikx]c[1] + x}^2\texttt{BesselY[2, -ikx]c[2]} \nonumber\\
	&\texttt{+}&\texttt{i/8(k}\pi\texttt{x}^3\texttt{BesselJ[0, ikx]BesselJ[2, ikx]BesselY[1, -ikx]} \nonumber \\
	&\texttt{+}& \texttt{k}\pi\texttt{x}^3\texttt{BesselJ[0, ikx]BesselJ[1, ikx]BesselY[2, -ikx])c[3]} \nonumber\\
	&\texttt{+}&\texttt{i/8(k}\pi\texttt{x}^3\texttt{BesselJ[2, ikx]BesselY[0, -ikx]BesselY[1, -ikx]} \nonumber\\
	&\texttt{+}& \texttt{k}\pi\texttt{x}^3\texttt{BesselJ[1, ikx]BesselY[0, -ikx]BesselY[2, -ikx])c[4]},\nonumber
	\end{eqnarray}
\end{subequations}
}which is not pleasant. Applying the brute force simplification algorithm gives
{\small
\begin{subequations} \label{ode:4_FullSimplify}
	\begin{eqnarray}
	\texttt{sol2} &=& \texttt{FullSimplify[sol]} \nonumber \\
	\texttt{y[x]} &\rightarrow& \texttt{-x}^2\texttt{/4(4BesselI[2, kx]c[1] - 4BesselY[2, -ikx]c[2]}\nonumber\\
	&\texttt{+}& \texttt{BesselI[0, kx]c[3] + BesselY[0, -ikx]c[4])} \nonumber,
	\end{eqnarray}
\end{subequations}}Using the relationship $Y_p(iz) \propto K_p(z)$ between Bessel and modified Bessel functions as well as the recurrence relations $I_0(z)=3I_1(z)/z + I_2(z)$ and $K_0(z)=3K_1(z)/z + K_2(z)$ converts \textit{Mathematica}'s simplified solution into our solution (\ref{ode:4_solution}).

\subsection{Repeated roots of Bessel's equation} \label{example5}

So far, we have solved problems that \textit{Mathematica} could manage, albeit more clumsily. Let us now tackle a problem that \textit{Mathematica} fails to solve directly, namely
\begin{equation} \label{ode:5_bessel}
\left(\frac{\mathrm{d}^2}{\mathrm{d}x^2} + \frac{1}{x}\frac{\mathrm{d}}{\mathrm{d}x} + 1\right)^3y = 0.
\end{equation}
Sixth-order equations do have applications, for example when studying the interactions between fluids and elastic media \cite{duprat2011dynamics}. In this case, we construct the homogeneous solutions by solving the sub-problems $\left(\frac{\mathrm{d}^2}{\mathrm{d}x^2} + \frac{1}{x}\frac{\mathrm{d}}{\mathrm{d}x} + 1\right)^2u_1=0$ and $\Big(\frac{\mathrm{d}^2}{\mathrm{d}x^2} + \frac{1}{x}\frac{\mathrm{d}}{\mathrm{d}x} + 1\Big)u_2=u_1$. To solve the first sub-problem for $u_1(x)$, we break things down further as before and solve the problems $\left(\frac{\mathrm{d}^2}{\mathrm{d}x^2} + \frac{1}{x}\frac{\mathrm{d}}{\mathrm{d}x} + 1\right)v_1 = 0$ and $\left(\frac{\mathrm{d}^2}{\mathrm{d}x^2} + \frac{1}{x}\frac{\mathrm{d}}{\mathrm{d}x} + 1\right)v_2 = v_1$ so that $u_1(x)=\mathrm{span}\{v_1(x),v_2(x)\}$.

The solutions for $v_1(x)$ are the Bessel functions $v_1(x)=c_1J_0(x) + c_2Y_0(x)$. Anticipating resonance in the equation for $v_2(x)$, we introduce $\epsilon$ and write
\begin{equation} \label{ode:5_2nd_subsub_problem}
\left(\frac{\mathrm{d}^2}{\mathrm{d}x^2} + \frac{1}{x}\frac{\mathrm{d}}{\mathrm{d}x} + 1\right)v_2 = c_1J_0\left((1+\epsilon)x\right) + c_2Y_0\left((1+\epsilon)x\right).
\end{equation}
Positing the particular solution $v_p = AJ_0\left((1+\epsilon)x\right) + BY_0\left((1+\epsilon)x\right)$ works only if we choose $A=\frac{c_1}{1-(1+\epsilon)^2}$ and $B=\frac{c_2}{1-(1+\epsilon)^2}$. Therefore, at this point the solution for $u_1(x)$ is
\begin{equation} \label{ode:5_u1_epsilon}
u_1(x) = c_1J_0(x) + c_2Y_0(x) + c_3\frac{J_0\left((1+\epsilon)x\right)}{1-(1+\epsilon)^2} + c_4\frac{Y_0\left((1+\epsilon)x\right)}{1-(1+\epsilon)^2},
\end{equation}
where we have relabelled some constants. Taylor expanding $J_0(x+\epsilon x) = J_0(x) - \epsilon xJ_1(x) + \mathrm{O}\left(\epsilon^2\right)$ and $Y_0(x+\epsilon x) = Y_0(x) - \epsilon xY_1(x) + \mathrm{O}\left(\epsilon^2\right)$, substituting them into equation (\ref{ode:5_u1_epsilon}), and keeping terms only up to first order in $\epsilon$ gives
\begin{equation} \label{ode:5_u1_continued}
u_1(x) = \left(c_1 - \frac{c_3}{2\epsilon}\right)J_0(x) + \left(c_2 - \frac{c_4}{2\epsilon}\right)Y_0(x) + \frac{c_3}{2}xJ_1(x) + \frac{c_4}{2}xY_1(x).
\end{equation}
Relabelling constants as usual gives the final solution for $u_1$
\begin{equation} \label{ode:5_u1_final}
u_1(x) = \tilde{c}_1J_0(x) + \tilde{c}_2Y_0(x) + \tilde{c}_3xJ_1(x) + \tilde{c}_4xY_1(x).
\end{equation}

We can now proceed to solve the sub-problem for $u_2(x)$. We have
\begin{equation} \label{ode:5_u2}
\left(\frac{\mathrm{d}^2}{\mathrm{d}x^2} + \frac{1}{x}\frac{\mathrm{d}}{\mathrm{d}x} + 1\right)u_2 = \tilde{c}_1J_0(x) + \tilde{c}_2Y_0(x) + \tilde{c}_3xJ_1(x) + \tilde{c}_4xY_1(x).
\end{equation}
As usual, we proceed by inserting the parameter $\epsilon$ such that
\begin{subequations} \label{ode:5_u2_continued}
\begin{eqnarray}
\left(\frac{\mathrm{d}^2}{\mathrm{d}x^2} + \frac{1}{x}\frac{\mathrm{d}}{\mathrm{d}x} + 1\right)u_2 = \tilde{c}_1J_0((1+\epsilon)x) + \tilde{c}_2Y_0((1+\epsilon)x) \\ + \tilde{c}_3xJ_1((1+\epsilon)x) + \tilde{c}_4xY_1((1+\epsilon)x) \nonumber
\end{eqnarray}
\end{subequations}
and propose a particular solution $u_p = AJ_0((1+\epsilon)x) + BY_0((1+\epsilon)x) + CxJ_1((1+\epsilon)x) + DxY_1\left((1+\epsilon)x\right)$. Substituting this anzats, we find after some simplification and rearrangement that
\begin{subequations} \label{ode:5_coefficients}
	\begin{eqnarray}
	A &=& \frac{\tilde{c}_1-2(1+\epsilon)C}{1-(1+\epsilon)^2} \\
	B &=& \frac{\tilde{c}_2-2(1+\epsilon)D}{1-(1+\epsilon)^2} \\
	C &=& \frac{\tilde{c}_3}{1-(1+\epsilon)^2} \\
	D &=& \frac{\tilde{c}_4}{1-(1+\epsilon)^2}.
	\end{eqnarray}
\end{subequations}
The full solution to the original problem (\ref{ode:5_bessel}) is $y=\mathrm{span}\{u_1,u_2\}$. Combining our results from equations (\ref{ode:5_u1_final}) and (\ref{ode:5_coefficients}) gives
\begin{subequations} \label{ode:5_continued_solution}
	\begin{eqnarray}
	~~~~~~~~~~y(x) &=& \tilde{c}_1J_0(x) + \tilde{c}_2Y_0(x) + \tilde{c}_3xJ_1(x) + \tilde{c}_4xY_1(x) \\ &+& AJ_0((1+\epsilon)x) + BY_0((1+\epsilon)x)
	+ c_5\frac{xJ_1((1+\epsilon)x)}{1-(1+\epsilon)^2} + c_6\frac{xY_1((1+\epsilon)x)}{1-(1+\epsilon)^2}.\nonumber
	\end{eqnarray}
\end{subequations}
The simplest way to proceed is to set $\epsilon=0$ in the terms involving $A$ and $B$, because we have freedom to make the transformations $\tilde{c}_1=\hat{c}_1 - A$ and $\tilde{c}_2=\hat{c}_2 - B$ to remove divergence issues. Thus we have
\begin{equation} \label{ode:5_continued_solution_2}
y(x) = \hat{c}_1J_0(x) + \hat{c}_2Y_0(x) + \tilde{c}_3xJ_1(x) + \tilde{c}_4xY_1(x) - c_5\frac{xJ_1((1+\epsilon)x)}{2\epsilon + \epsilon^2} - c_6\frac{xY_1((1+\epsilon)x)}{2\epsilon + \epsilon^2}.
\end{equation}
We can now perform the Taylor expansions
\begin{subequations}\label{ode:5_TaylorSeriesBesselFunctions}
\begin{eqnarray}
J_1(x + \epsilon x) &=& J_1(x) + \epsilon\left(J_1(x) - xJ_2(x)\right) + \mathrm{O}\left(\epsilon^2\right)\\
\hbox{and}~~Y_1(x + \epsilon x) &=& Y_1(x) + \epsilon\left(Y_1(x) - xY_2(x)\right) + \mathrm{O}\left(\epsilon^2\right),
\end{eqnarray}
\end{subequations}
where we have used the formulas $J_1'(z) = J_1(z)/z - J_2(z)$ and $Y_1'(z) = Y_1(z)/z - Y_2(z)$. Substituting these expansions in equation (\ref{ode:5_continued_solution_2}) while retaining terms only up to linear order in $\epsilon$ furnishes
\begin{subequations} \label{ode:5_continued_solution_3}
\begin{eqnarray}
~~~~~~~~~~y(x) &=& \hat{c}_1J_0(x) + \hat{c}_2Y_0(x) + \left(\tilde{c}_3-\frac{c_5}{2\epsilon} - \frac{c_5}{2}\right)xJ_1(x) + \left(\tilde{c}_4-\frac{c_6}{2\epsilon} - \frac{c_6}{2}\right)xY_1(x) \\ &+& \frac{c_5}{2}x^2J_2(x) + \frac{c_6}{2}x^2Y_2(x). \nonumber
\end{eqnarray}
\end{subequations}
Relabelling constants to clean up the solution and remove divergences finally gives all six linearly independent solutions to equation (\ref{ode:5_bessel}):
\begin{equation} \label{ode:5_final_solution}
y(x) = \hat{c}_1J_0(x) + \hat{c}_2Y_0(x) + \hat{c}_3xJ_1(x) + \hat{c}_4xY_1(x) + \tilde{c}_5x^2J_2(x) + \tilde{c}_6x^2Y_2(x).
\end{equation}

While we did all this work to show the explicit steps, the general structure of our method should now be clear. Problems of the form $\hat{D}^n[y(x)]=0$ can be broken down into $n-1$ resonance problems, the solutions to which can be generated from the function $u(x)$ where $\hat{D}[u(x)]=0$ by constructing the Taylor expansion 
\begin{equation} \label{ode:5_general_expansion}
u(x;\epsilon)=u(x) + \frac{\partial u}{\partial\epsilon}\bigg{|}_{\epsilon=0}\epsilon + \frac{\partial^2 u}{\partial\epsilon^2}\bigg{|}_{\epsilon=0}\frac{\epsilon^2}{2}+\cdots+\frac{\partial^{n-1}u}{\partial\epsilon^{n-1}}\bigg{|}_{\epsilon=0}\frac{\epsilon^{n-1}}{(n-1)!}.
\end{equation}
Therefore the function $\frac{\partial^k u}{\partial\epsilon^k}\big{|}_{\epsilon=0}$ gives the $k^{\mathrm{th}}$ repeated root of $\hat{D}^n[y(x)]=0$. In this example, $u(x;\epsilon) = J_0\left((1+\epsilon)x\right) + Y_0\left((1+\epsilon)x\right)$, and applying the formula (\ref{ode:5_general_expansion}) furnishes the delightful formula 
\begin{equation} \label{ode:5_kth_root_Bessel}
    u_{(k)}(x) = b_{1}x^kJ_k(x) + b_{2}x^kY_k(x)
\end{equation}
for the $k^{\mathrm{th}}$ repeated root of Bessel's equation.

\subsection{Legendre's equation at resonance} \label{example6}

We now turn to an example where the manner in which we introduce the parameter $\epsilon$ is quite different. We consider the following Legendre BVP at resonance
\begin{subequations} \label{ode:6-legendre}
	\begin{eqnarray}
	\frac{\mathrm{d}}{\mathrm{d}x}\left[\left(1-x^2\right)\frac{\mathrm{d}y}{\mathrm{d}x}\right] + n(n+1)y &=& P_n(x) \\
	y(1) &=& 1 \\
	y(x) &=&\mathrm{finite~on}~x\in(-1, 1],
	\end{eqnarray}
\end{subequations}
where $P_n(x)$ is the Legendre polynomial of the first kind of order $n$ for $n\in\mathbb{Z}$. In this case, the homogeneous solutions are $u_1(x) = P_n(x)$ and $u_2(x) = Q_n(x)$, which are the Legendre functions of the first and second kind, respectively. 

Clearly, we have resonance, but if we try our usual strategy of perturbing the forcing function as $P_n\left((1+\epsilon)x\right)$ and utilize the ansatz $u_p(x)=AP_n\left((1+\epsilon)x\right)$, it will fail because of the structure of the $(1-x^2)$ term in equation (\hyperref[ode:6-legendre]{3.44a}). Instead, for equations of this type it is more useful to perturb the \textit{order} of the forcing function, in this case the order of the Legendre polynomial $n$. We have
\begin{equation} \label{ode:6-legendre_perturbed}
    \frac{\mathrm{d}}{\mathrm{d}x}\left[\left(1-x^2\right)\frac{\mathrm{d}y}{\mathrm{d}x}\right] + n(n+1)y = P_{n+\epsilon}(x),
\end{equation}
where we interpret $P_{n+\epsilon}(x)$ as a homogeneous solution to equation (\hyperref[ode:6-legendre]{3.44a}) for $n\rightarrow n+\epsilon\in \mathbb{R}$, which can be represented by a hypergeometric function \cite[eq. 15.9.7]{NIST:DLMF}, although we will not make use of this directly. Positing the ansatz $u_p(x) = AP_{n+\epsilon}(x)$ works and we find $A=\frac{1}{n(n+1) - (n+\epsilon)(n+\epsilon+1)}$, and so the general solution is
\begin{equation} \label{ode:6-legendre_perturbed_general_solution}
    y(x) = c_1P_n(x) + c_2Q_n(x) - \frac{P_{n+\epsilon}(x)}{(2n+1)\epsilon+\epsilon^2}.
\end{equation}
We now Taylor expand the particular solution in the usual way, except now we are expanding in the order of the Legendre polynomial. Keeping terms only up to $\mathrm{O}(\epsilon)$ we obtain
\begin{equation} \label{ode:6-taylor_series}
    -u_p(x) = \frac{P_n(x)}{(2n+1)\epsilon} + \frac{1}{2n+1}\frac{\partial P_n(x)}{\partial n}.
\end{equation}
We now have to evaluate the derivative $\frac{\partial P_n(x)}{\partial n}$, which we do by thinking of $n$ as a continuous variable, then evaluating the result only at integer $n$. There are many ways to do this, for example in terms of the aforementioned hypergeometric function, or by using an integral representation valid for arbitrary $n$ \cite[eq. 4.1]{szmytkowski2006derivative}. There is actually a simpler way of computing $\frac{\partial P_n(x)}{\partial n}$, which was discovered by Jolliffe \cite{jolliffe1919form}. His method leads to the remarkable formula
\begin{equation} \label{ode:6-legendre_derivative}
    \frac{\partial P_n(x)}{\partial n} := P_{n,1}(x) =  \frac{1}{2^{n-1}n!}\frac{\mathrm{d}^n}{\mathrm{d}x^n}\left[(x^2-1)^n\log\left(\frac{x+1}{2}\right)\right] - P_n(x)\log\left(\frac{x+1}{2}\right),
\end{equation}
which is interesting in that it is almost the standard Rodrigues formula $P_n(x) = \frac{1}{2^{n}n!}\frac{\mathrm{d}^n}{\mathrm{d}x^n}(x^2-1)^n$ apart from logarithmic corrections. Equipped with equation (\ref{ode:6-legendre_derivative}), the general solution becomes
\begin{equation} \label{ode:6_perturbed_general_solution}
    y(x) = \left(c_1-\frac{1}{(2n+1)\epsilon}\right)P_n(x) + c_2Q_n(x) - \frac{P_{n,1}(x)}{2n+1}.
\end{equation}
Equation (\ref{ode:6_perturbed_general_solution}) has the same structure as the solution derived by Backhouse \cite[eqs. 13-14]{backhouse2001resonant} in his study of the resonant Legendre equation, where he too remarked on the appearance of logarithmic terms.

We can now apply the boundary data. Condition (\hyperref[ode:6-legendre]{3.44c}) forces $c_2=0$ while condition (\hyperref[ode:6-legendre]{3.44b}) gives $1 = \left(c_1-\frac{1}{(2n+1)\epsilon}\right)P_n(1) - \frac{P_{n,1}(1)}{2n+1}$. Using the fact that $P_n(1)=1$, one can see that $P_{n,1}(1)=0$ due to the logarithm terms, hence $c_1 = 1 + \frac{1}{(2n+1)\epsilon}$ and the solution to the BVP (\hyperref[ode:6-legendre]{3.44}) is
\begin{equation} \label{ode:6-final_solution}
    y(x) = P_n(x) - \frac{P_{n,1}(x)}{2n+1},
\end{equation}
where $P_{n,1}(x)$ is given by equation (\ref{ode:6-legendre_derivative}).

\subsection{Repeated roots of the Hermite equation} \label{example7}

As our final example we now consider repeated roots of Hermite's equation
\begin{equation} \label{ode:7-hermite} 
    \left(\frac{\mathrm{d}^2}{\mathrm{d}x^2}-2x\frac{\mathrm{d}}{\mathrm{d}x} + 2n\right)^2y=0
\end{equation}
for $n\in\mathbb{Z}$. As usual, we seek to solve the sub-problems $\left(\frac{\mathrm{d}^2}{\mathrm{d}x^2}-2x\frac{\mathrm{d}}{\mathrm{d}x} + 2n\right)u_1 = 0$ and $\Big(\frac{\mathrm{d}^2}{\mathrm{d}x^2}-2x\frac{\mathrm{d}}{\mathrm{d}x} + 2n\Big)u_2 = u_1$. The solutions to the first sub-problem are the Hermite functions $u_1(x) = c_1H_n(x) + c_2G_n(x)$, where $H_n(x)$ is the Hermite polynomial of order $n$ and $G_n(x)=H_n(x)\int^{x}\mathrm{d}x'e^{x'^2}/H_n^2(x')$\footnote{There is no ``Hermite function of the $2^{\mathrm{nd}}$ kind'' commonly cited in the literature, so we give an integral representation here that can be derived using reduction of order (see \hyperref[Apdx-2]{Appendix 5.2}).}. Therefore, the second sub-problem we need to solve is 
\begin{equation} \label{ode:7-subproblem} 
    \left(\frac{\mathrm{d}^2}{\mathrm{d}x^2}-2x\frac{\mathrm{d}}{\mathrm{d}x} + 2n\right)u_2 = c_1H_n(x) + c_2G_n(x),
\end{equation}
which is a resonance problem. As in the previous example, we perturb the order of the Hermite functions and construct the perturbed problem
\begin{equation} \label{ode:7-perturbed} 
    \left(\frac{\mathrm{d}^2}{\mathrm{d}x^2}-2x\frac{\mathrm{d}}{\mathrm{d}x} + 2n\right)u_2 = c_1H_{n+\epsilon}(x) + c_2G_{n+\epsilon}(x),
\end{equation}
where we interpret $H_{n+\epsilon}(x)$ and $G_{n+\epsilon}(x)$ as a homogeneous solutions to equation (\ref{ode:7-hermite}) for $n+\epsilon\in\mathbb{R}$. As usual, we propose the form $u_p(x) = AH_{n+\epsilon}(x) + BG_{n+\epsilon}(x)$ and substitute into equation (\ref{ode:7-perturbed}). We find this works as long as we choose $A=-c_1/2\epsilon$ and $B=-c_2/2\epsilon$. Hence the general solution so far is 
\begin{equation} \label{ode:7-perturbed_solution} 
    y(x) = c_1H_n(x) + c_2G_n(x) + c_3\frac{H_{n+\epsilon}(x)}{\epsilon} + c_4\frac{G_{n+\epsilon}(x)}{\epsilon},
\end{equation}
where we have relabeled some integration constants. To finish, we construct the Taylor expansions $H_{n+\epsilon}(x) = H_{n}(x) + \frac{\partial H_n(x)}{\partial n}\epsilon + \mathrm{O}(\epsilon^2)$ and $G_{n+\epsilon}(x) = G_{n}(x) + \frac{\partial G_n(x)}{\partial n}\epsilon + \mathrm{O}(\epsilon^2)$, where the task now becomes computing $\frac{\partial H_n(x)}{\partial n}$, from which $\frac{\partial G_n(x)}{\partial n}$ follows by invoking the chain rule. Unfortunately this derivative appears less often in practice than the Legendre case, so there are no simple formulas commonly used. Thus we give an answer in terms of the confluent hypergeometric function $_1F_1$ \cite{WeissteinHermite}, 
\begin{equation} \label{ode:7-derivative_hermite_1}
    \frac{\partial H_n(x)}{\partial n} := H_{n,1}(x) = x\frac{\mathrm{d}}{\mathrm{d}\xi}\left[_1F_1\left(\xi, \frac{3}{2}, x^2\right)\right]\bigg{|}_{\xi=\frac{1-n}{2}}.
\end{equation}
Using this definition and the chain rule we also find
\begin{equation} \label{ode:7-derivative_hermite_2}
    \frac{\partial G_n(x)}{\partial n} := G_{n,1}(x) = H_{n,1}(x)\int^{x}\mathrm{d}x'~\frac{e^{x'^2}}{H_n^2(x')} - 2H_n(x)\int^{x}\mathrm{d}x'~\frac{e^{x'^2}H_{n,1}(x')}{H_n^3(x')}.
\end{equation}
With everything known, we substitute these functions into their respective Taylor expansions, group terms of order $\epsilon^{-1}$ with the homogeneous solutions, relabel integration constants, and set $\epsilon=0$. The final result is
\begin{equation} \label{ode:7-solution}
    y(x) = \tilde{c}_1H_n(x) + \tilde{c}_2G_n(x) + c_3H_{n,1}(x) + c_4G_{n,1}(x).
\end{equation}

\section{General structure and final remarks}

In this paper we highlighted how to construct resonant or repeated root solutions to ODEs via analytic continuation of an already known homogeneous solution in an introduced parameter $\epsilon$. For the practitioner who comes across these problems in their work, we hope these examples have offered a concrete guide on how to find such solutions. To conclude, we now give a more general derivation that encompasses all our examples.

We are interested in finding a particular solution $u_p(x)$ of the resonant ODE $\hat{L}\left[u_p(x)\right]=u(x)$, where $\hat{L}\left[u(x)\right]=0$. Suppose that $\hat{L}$ can be written as $\hat{L} = \hat{M} - \lambda$, where $\lambda$ is either a parameter in the original problem, or it can be introduced to the problem and set to unity at the end. The homogeneous solution will therefore depend on this parameter, and we denote this dependence by $u=u(x;\lambda)$. Therefore we seek to solve
\begin{equation} \label{general-1}
    \left(\hat{M}-\lambda\right)\left[u_p(x)\right] = u(x;\lambda).
\end{equation}
We now perturb this problem by introducing a small parameter $\epsilon$
\begin{equation} \label{general-2}
    \left(\hat{M}-\lambda\right)\left[u_p(x)\right] = u(x;\lambda + \epsilon)
\end{equation}
and propose that the particular solution has the form $u_p(x) = Au(x; \lambda + \epsilon)$. Since $\hat{M}[u(x;\lambda + \epsilon)]=(\lambda + \epsilon)u(x;\lambda + \epsilon)$, this ansatz satisfies equation (\ref{general-2}) if we choose $A = 1/\epsilon$, hence $u_p(x)=u(x;\lambda + \epsilon)/\epsilon$. We now Taylor expand $u_p(x)$ around $\epsilon=0$ to obtain\footnote{Note that $\frac{\partial u}{\partial\lambda} = \frac{\partial u}{\partial\epsilon}\big{|}_{\epsilon=0}$ identically, so this is consistent with our previous notation up until now.}
\begin{equation} \label{general-3}
    u_p(x) = \frac{u(x)}{\epsilon} + \frac{\partial u}{\partial\lambda} + \mathrm{O}\left(\epsilon\right). 
\end{equation}
The $u(x)/\epsilon$ term is proportional to the homogeneous solution $u(x)$, so we are free to remove it by lumping it with the integration constant associated with $u(x)$. This removes the $1/\epsilon$ divergence and allows us to set $\epsilon = 0$, resulting in
\begin{equation} \label{general-4}
    u_p(x) = \frac{\partial u}{\partial\lambda}.
\end{equation}
With this general derivation, we see that the resonant solution to an ODE is \textit{always} given by the derivative with respect to the eigenvalue $\lambda$ of the governing differential operator $\hat{L}$. We note that equation (\ref{general-4}) has also been derived by Makarov et al. \cite{makarov1978construction}. See \hyperref[table1]{Table 4.1} for a guide on how to apply this equation for the differential operators discussed in our examples.

\begin{table}[h!]
\centering
\caption{Construction of the resonant particular solution $u_p(x)$ from its associated homogeneous solution $u(x)$ using equation ($\ref{general-4}$). The calculation is done using the chain rule $\frac{\partial u}{\partial\lambda} = \frac{\mathrm{d}\mu}{\mathrm{d}\lambda}\frac{\partial u}{\partial\mu}$. To generate higher order repeated root solutions, simply take higher derivatives with respect to $\lambda$ as per equation ($\ref{general-5}$).}
\begin{tabular}{| P{11.0em} | P{4.1em} | P{6.5em} | P{9.0em} |} 
 \hline
 $\hat{L}$ & $\lambda$ & Homogeneous Solution $u(x)$ & Resonant Solution $u_p(x)$ \\ [1ex] 
 \hline\hline
 $\frac{\mathrm{d}^2}{\mathrm{d}x^2} + \mu^2$ & $-\mu^2$ & $\sin\mu x$, $\cos\mu x$ & $-\frac{x\cos\mu x}{2\mu}$, $\frac{x\sin\mu x}{2\mu}$ \\ [1ex] 
 \hline
 $x\frac{\mathrm{d}}{\mathrm{d}x} - \mu$ & $\mu$ & $x^\mu$ & $x^{\mu} \log x$ \\ [1ex] 
 \hline
 $\frac{1}{x}\frac{\mathrm{d}^2}{\mathrm{d}x^2} - \mu^3$ & $\mu^3$ & $\mathrm{Ai}(\mu x)$, $\mathrm{Bi}(\mu x)$ & $\frac{x\mathrm{Ai}'(\mu x)}{3\mu^2}$, $\frac{x\mathrm{Bi}'(\mu x)}{3\mu^2}$ \\ [1ex] 
 \hline
 $\frac{\mathrm{d}^2}{\mathrm{d}x^2} + \frac{1}{x}\frac{\mathrm{d}}{\mathrm{d}x} + \mu^2$ & $-\mu^2$ & $J_0(\mu x)$, $Y_0(\mu x)$ & $\frac{xJ_1(\mu x)}{2\mu}$, $\frac{xY_1(\mu x)}{2\mu}$ \\ [1ex] 
 \hline
 $\frac{\mathrm{d}}{\mathrm{d}x}\left[(1-x^2)\frac{\mathrm{d}}{\mathrm{d}x}\right] + \mu(\mu+1)$ & $-\mu(\mu+1)$ & $P_{\mu}(x)$, $Q_{\mu}(x)$ & $-\frac{P_{\mu, 1}(x)}{2\mu + 1}$, $-\frac{Q_{\mu, 1}(x)}{2\mu + 1}$ \\ [1ex] 
 \hline
 $\frac{\mathrm{d}^2}{\mathrm{d}x^2} - 2x\frac{\mathrm{d}}{\mathrm{d}x} + 2\mu$ & $-2\mu$ & $H_{\mu}(x)$, $G_{\mu}(x)$ & $-\frac{H_{\mu, 1}(x)}{2}$, $-\frac{G_{\mu, 1}(x)}{2}$ \\ [1ex] 
 \hline
\end{tabular}
\label{table1}
\end{table}

From example \hyperref[example5]{3.5}, we showed that the repeated root problem $\hat{D}^k[y(x)] = 0$ is equivalent to $k - 1$ resonance problems. Thus in the context of this derivation, the $k^{\mathrm{th}}$ repeated root solution is given by 
\begin{equation} \label{general-5}
    u_{(k)}(x) = \frac{\partial^k u}{\partial\lambda^k}
\end{equation}
where $\hat{D}[u(x)]=0$.

While all the problems presented here may be solved using reduction of order, the benefit of using the method presented here should be apparent. There is minimal algebra and no integration of a reduced-order ODE is  required. It is only necessary to compute a Taylor series. In this way, the connection between the homogeneous solutions and resonant or repeated root solutions becomes clear. We hope that this note motivates instructors to present this approach along side reduction of order in their differential equation courses. These problems do arise in practice, and we believe this method is helpful in exercising thought processes useful in applied mathematics, as well as producing the most elegant solution. 

\pagebreak

\section{Appendix} 

\subsection{Calculation of $\mathrm{Ai}(0)$} \label{Apdx-1}

\begin{figure}[H] \label{fig:1}
\setlength{\belowcaptionskip}{-10pt}
\centering
\includegraphics[width=3.5in, height=2.5in]{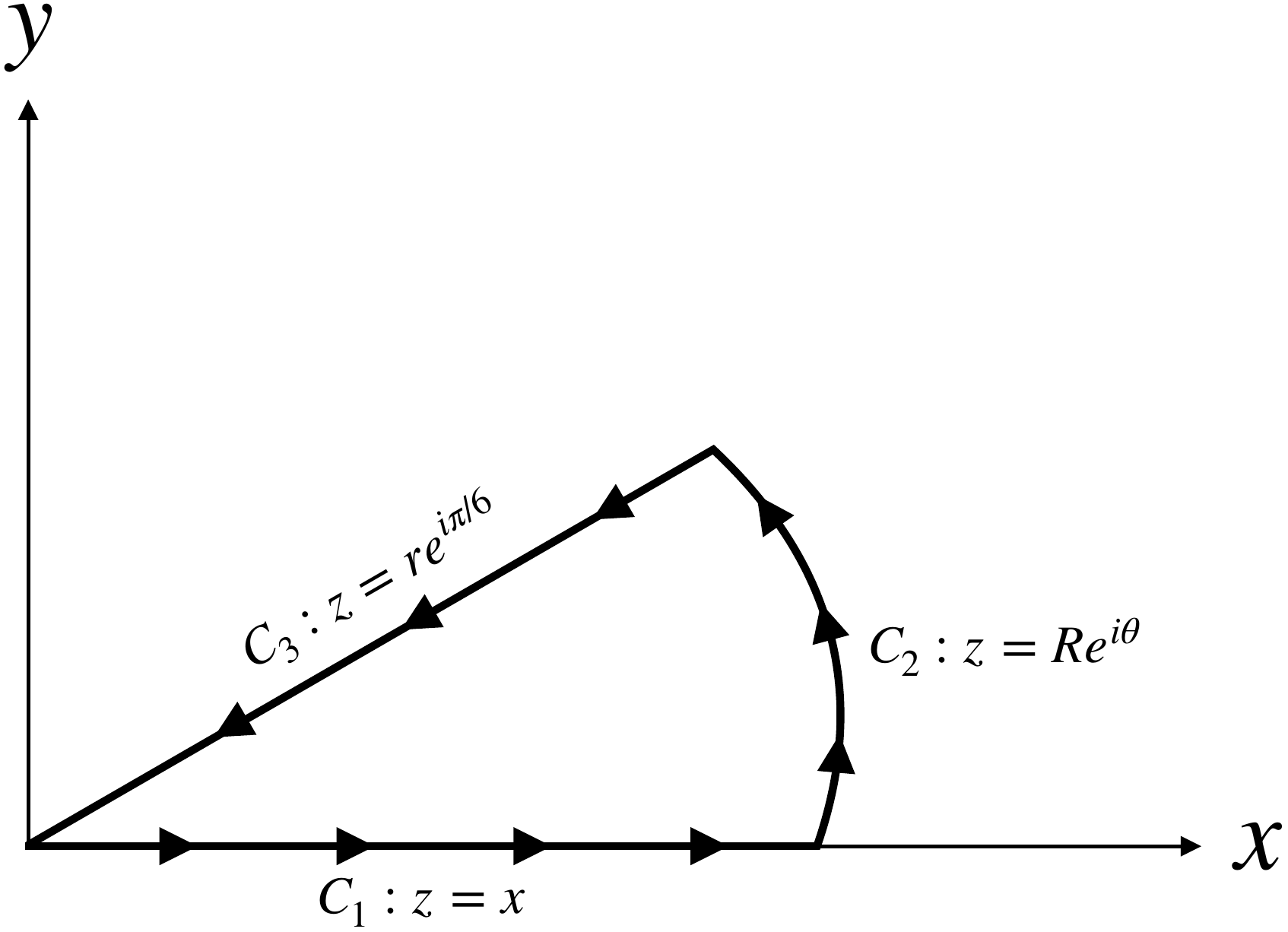}
\caption{Contour used for evaluation of $I=\oint_{C}\mathrm{d}z~e^{iz^3}$. $C_3$ is chosen as the direction of steepest descent for the integrand.}
\end{figure}

We start with the well-known integral representation $\mathrm{Ai}(x)=\frac{1}{\pi}\int_0^\infty\mathrm{d}t\cos\left(\frac{t^3}{3} + xt\right)$. Our method of approach in evaluating $\mathrm{Ai}(0)$ will be standard contour integration, so instead we consider the complex integral
\begin{equation} \label{contour-1}
    \mathcal{I}(\alpha) := \int_0^{\infty}\mathrm{d}w~e^{i\alpha w^3} = \frac{1}{\alpha^{1/3}}\underbrace{\int_0^{\infty}\mathrm{d}z~e^{i z^3}}_{:=I}.
\end{equation}
One can see that $\mathrm{Re}\left[\mathcal{I}(1/3)\right]=\pi\mathrm{Ai}(0)$, and so our task becomes evaluating $I$. Because $e^{iz^3}$ is entire, $\oint_{C}\mathrm{d}z~e^{iz^3} = 0$ by Cauchy's theorem. We choose the contour $C$ as show in \hyperref[fig:1]{Fig. 5.1}, where we will take $R\rightarrow\infty$. Writing out the integrals we have
\begin{equation} \label{contour-2}
    \lim_{R\rightarrow\infty}\left[\underbrace{\int_0^{R}\mathrm{d}x~e^{ix^3}}_{I} + \underbrace{\int_0^{\pi/6}Re^{i\theta}i\mathrm{d}\theta~e^{iR^3e^{i3\theta}}}_{J} + \underbrace{\int_R^{0}e^{i\pi/6}\mathrm{d}r~e^{-r^3}}_{K}\right] = 0,
\end{equation}
where $I$ is the integral we wish to compute. 

We can prove that $J$ vanishes in the $R\rightarrow\infty$ limit by the following argument:
\begin{subequations} \label{contour-3}
	\begin{eqnarray}
	|J| &=& \int_0^{\pi/6}\big{|}Re^{i\theta}i\mathrm{d}\theta~e^{iR^3e^{i3\theta}}\big{|} \\
	&\leq& R\int_0^{\pi/6}\big{|}\mathrm{d}\theta~e^{-R^3\sin 3\theta}\big{|} \\
	&\leq& R\int_0^{\pi/6}\mathrm{d}\theta~e^{-R^3 6\theta/\pi} \\
	&=& \frac{\pi}{6R^2}\left(1 - e^{-R^3}\right)\rightarrow0~~~\mathrm{as}~~~R\rightarrow\infty.
	\end{eqnarray}
\end{subequations}
In the first step, we get rid of all pure unimodular phases. In the next step, we use the bound $\sin 3\theta>6\theta/\pi$ on $\theta\in[0,\pi/6]$, which allows us to evaluate the integral and prove that it vanishes. Hence, $I = -K$, which we can evaluate in terms of the Gamma function as follows:
\begin{subequations} \label{contour-4}
	\begin{eqnarray}
	I &=& -K \\
	&=& e^{i\pi/6}\int_0^{\infty}\mathrm{d}r~e^{-r^3} \\
	&=& \frac{e^{i\pi/6}}{3}\int_0^{\infty}\mathrm{d}u~u^{-2/3}e^{-u} \\
	&=& \frac{e^{i\pi/6}}{3}\Gamma(1/3).
	\end{eqnarray}
\end{subequations}
Here we have just made the change of variables $u=r^3$ and utilized the definition of the Gamma function. 

Hence we have $\pi\mathrm{Ai}(0)=\mathrm{Re}[\mathcal{I}(1/3)]=3^{1/3}\mathrm{Re}[I]=\frac{3^{1/3}\cos(\pi/6)\Gamma(1/3)}{3}=\frac{3^{5/6}\Gamma(1/3)}{6}$. Simplifying a bit gives the final result $\mathrm{Ai}(0)=\frac{\Gamma(1/3)}{2\pi 3^{1/6}}$.

\subsection{Calculation of the Hermite function of the $2^{\mathrm{nd}}$ kind $G_n(x)$} \label{Apdx-2}

Let $v_2(x)=H_n(x)f(x)$ be the $2^{\mathrm{nd}}$ linearly independent homogeneous solution to Hermite's equation $v_2''-2xv_2'+2nv_2=0$, where $f(x)$ is to be determined. Substituting this ansatz into Hermite's equation and simplifying a bit results in
\begin{equation} \label{roo-1}
    H_nf'' + 2H_n'f' - 2xf'H_n + \underbrace{\left(H_n'' - 2xH_n' + 2nH_n\right)}_{=0}f=0
\end{equation}
which is a first order ODE for $f'$. We can divide equation (\ref{roo-1}) through by $f'H_n$ and then rearrange it into the form
\begin{equation} \label{roo-2}
    \frac{\mathrm{d}}{\mathrm{d}x}\log f'H_n^2 = 2x,
\end{equation}
which can be integrated to give $f'(x)=e^{x^2}/H_n^2(x)$ (neglecting the integration constant). Another integration produces $f(x)=\int^{x}\mathrm{d}x'e^{x'^2}/H_n^2(x')$, and thus our desired result is
\begin{equation} \label{roo-3}
    v_2(x) = H_n(x)\int^{x}\mathrm{d}x'\frac{e^{x'^2}}{H_n^2(x')}.
\end{equation}

\section*{Acknowledgements}

B.G. is supported by the Paul and Daisy Soros Fellowship and the NSF Graduate Research Fellowship Program. H.A.S. recalls learning this approach for constant coefficient ODEs forced at resonance when first teaching an undergraduate course 30 years ago while following notes from John Hutchinson. We thank Peter Howell, Dionisios Margetis, Ali Nadim, Elie Raphael, and Tom Witelski for their helpful comments and suggestions. 

\bibliography{references.bib}
\bibliographystyle{siamplain}

\end{document}